# Generalized Small Cancellation Theory

Günther Huck
Northern Arizona
University

Stephan Rosebrock
Johann–Wolfgang Goethe
Universität

October 26, 1995

**Abstract**

We present four generalized small cancellation conditions for finite presentations and solve the word- and conjugacy problem in each case. Our conditions $W$ and $W^*$ contain the non-metric small cancellation cases C(6), C(4)T(4), C(3)T(6) (see [7]) but are considerably more general. $W$ also contains as a special case the small cancellation condition $W(6)$ of Juhasz [6]. If a finite presentation satisfies $W$ or $W^*$ then it has a quadratic isoperimetric inequality and therefore solvable word problem. For the class $W$ this was first observed by Gersten in [3] which also contains an idea of the proof. Our main result here is the proof of the conjugacy problem for the classes $W^*$ and $W$ which uses the geometry of non-positively curved piecewise Euclidean complexes developed by Bridson in [1].

The conditions $V$ and $V^*$ generalize the small cancellation conditions C(7), C(5)T(4), C(4)T(5), C(3)T(7). If a finite presentation satisfies the condition $V$ or $V^*$, then it has a linear isoperimetric inequality and hence the group is hyperbolic.

## 1 The Conditions $W^*$ and $W$

Let $P = <x_1, \ldots, x_n \mid R_1, \ldots, R_m>$ be a finite presentation of the group $G$. We always assume that each relator $R_i$ is cyclically reduced and no relator is the trivial word or a cyclic permutation of another relator or of the inverse of another relator. (By "relator" we always mean a defining relator of the presentation.) Let $F$ be the free group on the generators. $K_P$ denotes the standard 2-complex modeled on $P$. The *Whitehead graph* $W_P$ of $K_P$ is the boundary of a regular neighborhood of the only vertex of $K_P$. For each generator $x_i$ of $P$ it has two vertices $+x_i$ and $-x_i$ which correspond to the beginning and the end of the oriented loop labeled $x_i$ in $K_P$. The edges of $W_P$ are the *corners* of the 2-cells of the 2-complex. The *star graph* $S_P$ is the same as the Whitehead graph if no relator of $P$ is a proper power. If $P$ contains relators that are proper powers, the stargraph is obtained from the Whitehead graph by identifying edges of $W_P$ that correspond under the periodicity of a relator $R_i$. More precisely, if the relators $R_i$ of $P$ are of the form $s_i^{k_i}$ with $s_i$ not a proper power, then the star graph $S_P$ is the Whitehead graph of the presentation $< x_1, \ldots, x_n \mid s_1, \ldots, s_m >$ ($s_i$ is called the *root* of the relator $R_i$). Let $\pi: W_P \to S_P$ be the natural map from the Whitehead graph to the stargraph. We use the star graph in dealing with algebraic questions about presentations, such as the word problem and the conjugacy problem, because, algebraically, a rotation of a periodic relator by any multiple of its period is irrelevant. (For topological questions, e.g. when dealing with asphericity of 2-complexes, such rotations are quite significant and therefore the Whitehead graph is used instead.)

A *combinatorial map* is a cellular map that maps each open cell homeomorphically onto an open cell. A cell complex is said to be *combinatorial* if its attaching maps are combinatorial. Clearly, a 2- complex $K_P$ modeled on a presentation $P$ is combinatorial. A *diagram over* $P$ is a combinatorial map $f: M \to K_P$, where $M = S^2 - \cup_{i \leq \tau} \sigma_i^2$ is a combinatorial cell decomposition of a 2-sphere minus $\tau$ disjoint open 2-cells. If $\tau = 0$, $M$ is a sphere and we speak of a *spherical diagram*. If $\tau = 1$, we call $M$ a *disk diagram* over $P$ which is the same as a van Kampen diagram. If $\tau = 2$ we speak of an *annular diagram* over $P$.

For $\tau > 0$ let $\delta M$ denote the boudary of $M$, where we think of the boundary as a collection of boundary paths, namely the boundary paths of the 2-cells $\sigma_i$ that are deleted from $S^2$ to obtain $M$. The *length of the boundary* $l(\delta M)$ is the number of edges in the boundary paths, counting with multiplicity edges that are traversed twice in opposite direction. An oriented edge $e$ in M is labeled by the name $x_i^\epsilon$ of its image $f(e)$ in $K_P$ ($\epsilon \in \{-1, +1\}$), and hence every edge path in $M$ is labeled by a word $w$ in $F$.

A word $w$ in $F$ is trivial in $G$ if and only if there is a disk diagram over $P$ reading $w$ along its boundary. Two words $u$ and $v$ in $F$ represent conjugate elements in $G$ if and only if there is an annular diagram over $P$ such that its boundary paths, oriented parallely, read for suitable choices of start points the words $u$ and $v$. (A parallel orientation is an orientation of the two boundary paths such that they are freely homotopic in $M$) Both facts are well known. Since the proof of the relation between conjugate words and annular diagrams is not quite as standard it will be sketched in section 4.

The Whitehead graph $W(K)$ of an arbitrary combinatorial 2-complex $K$ is the boundary of the regular neighborhood of the 0-skeleton. $W(K)$ is realized as a 1-complex as follows: Its vertices are the intersections of the boundary of the regular neighborhood with the 1-cells of $K$, and its open edges are the components of the intersections with the open 2-cells of $K$ (i.e. the corners of the 2-cells of $K$). A combinatorial map between 2-complexes induces a combinatorial map between the corresponding Whitehead graphs; this applies in particular to maps $f: M \to K_P$ which define diagrams over $P$. For a diagram $M$ over $P$ the components of the Whitehead graph $W(M)$ are either 1-manifolds (i.e. circles or arcs) or isolated vertices. We say the diagram $M$ is *vertex reduced* if each component $w \subset W(M)$ maps under $\pi \circ f$ onto a path in the stargraph $S_P$ that does not pass an edge twice in different directions.

If a diagram $M$ is not vertex reduced, an elementary reduction move can be performed to $M$, as indicated in figure 1. In this figure, $D_1$ and $D_2$ are distinct 2-cells of $M$ that intersect at least in one vertex $Q$ and map to the same 2-cell in $K_P$ with opposite orientation, such that the oppositely oriented corners $\alpha$ and $\beta$ at $Q$ map to the same edge $\gamma$ in $S_P$. (This is equivalent to the oppositely oriented boundary paths $\nu_1$ and $\nu_2$ mapping to the same path labeled $w$ in $K_P$.) The elementary reduction deletes the 2-cells $D_1$ and $D_2$ in $M$ and identifies the paths $\nu_1$ and $\nu_2$ to the path $\nu$ as sketched in figure 1.

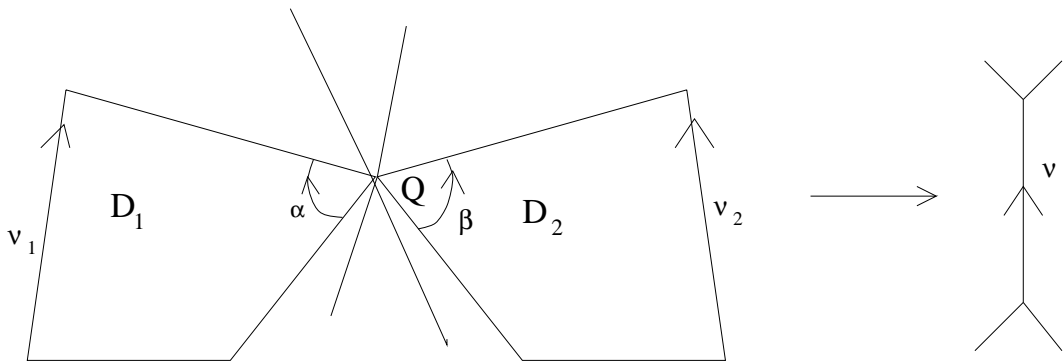

Figure 1: Reduction along a vertex

In general, the situation may be a little more complicated than figure 1 suggests. An elementary reduction where the boundaries of $D_1$ and $D_2$ intersect in more than one component or where the boundary of each $D_i$ has self intersections may cause parts of the diagram to be squeezed off (mostly as spheres). More precisely, let the boundary paths $\nu_1$ and $\nu_2$ be given by the sequences of vertices and oriented edges in $M: Q = Q_0 a_1 Q_1 a_2 Q_2 ... a_n Q_n = Q$ and $Q = R_0 b_1 R_1 b_2 ... b_n R_n = Q$, respectively. If $a_i = b_i$ or $Q_i = R_i$ we speak of a symmetric identification of edges or vertices, and if $(a_i = a_j$ and $b_i = b_j)$ or $(Q_i = Q_j$ and $R_i = R_j)$ we speak of a symmetric pair of self identifications. Let $c$ be the number of components of symmetric identifications and $c'$ the number of components of symmetric pairs of self identifications (identifications of connected segments are counted as one). Then the elementary reduction will have the effect of "squeezing off" $c + c'$ components from the diagram, i.e. the result are $c + c' + 1$ diagrams, joined by successively taking one-point-unions. In case $M$ is a disk diagram, one of these $c + c' + 1$ components will contain the boundary of $M$ (which as a path labeled by a word $w$ in $F$ remains unchanged), the others will be spherical. By deleting the spherical components we obtain a simpler diagram for the same boundary word. Now let $M$ be an annular diagram that is split by an elementary reduction into several components with respect to one-point-union. If the boundary paths which represent the conjugate words $u$ and $v$ both end up in the same component, we can delete the remaining spherical components and obtain a smaller annular diagram for $u$ and $v$. If the two boundary paths end up in different components (which are disk diagrams), $u$ and $v$ will both be equal to 1 in $G$. Since standardly a proof of the conjugacy problem is based on a solution of the word problem, this represents a trivial case.

This splitting off of components under an elementary reduction is easier to observe if one uses, instead of diagrams, the equivalent concept of pictures. A picture is essentially a homotopic representation of the map $f: M \to K_P$ by a transverse regular map (for details see [9]). We do not introduce this otherwise useful concept here since we do not need it for the proofs of the following theorems.

By iterating the above procedure of elementary reduction followed by deleting the spherical components, any disk diagram for a word $w$ and any annular diagram for conjugate words $u$ and $v$ (that are not equal to 1 in $G$) can be transformed into a corresponding vertex reduced diagram that has fewer 2-cells.

We also need the more general notion of a "reduced diagram" which is standardly used in small cancellation theory and other aplications of diagrams over groups. A diagram is said to be *reduced* if it does not contain a pair of (distinct) 2-cells $D_1$ and $D_2$ that intersect in at least one edge $x$ (see figure 2) such that their oppositely oriented boundary paths $x\nu_1^{-1}$ and $x\nu_2^{-1}$ read the same word in $F$. Equivalently, we can say $M$ is reduced if no component of the Whitehead graph of $M$ maps under $\pi \circ f$ onto a path in $S_P$ with backtracking, i.e. a path that travels an edge $\gamma$ immediately followed by its inverse. If a diagram $M$ is not reduced, an elementary reduction can be performed to $M$ that deletes the pair of 2-cells $D_1$ and $D_2$ and identifies the paths $\nu_1$ and $\nu_2$. Regarding the details, the same discussion as above applies.

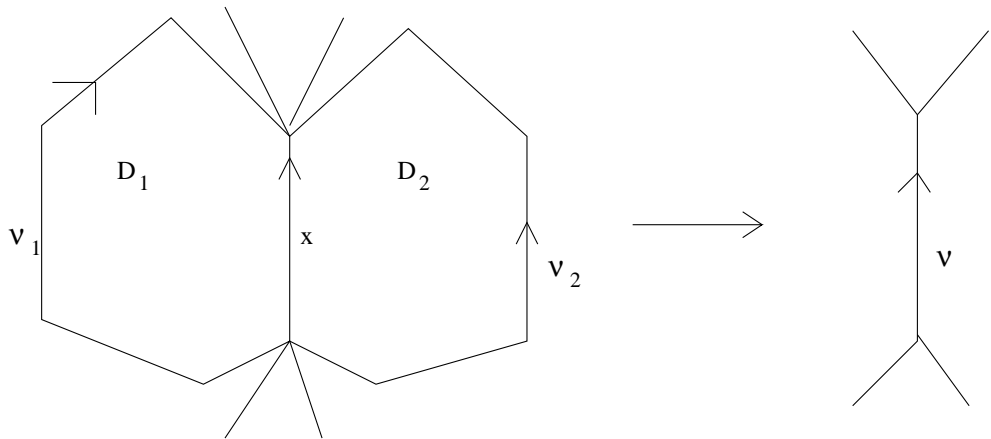

Figure 2: Reduction along an edge

the standard notion of "reduced diagram" will allow us to define the small cancellation classes $W^*$ and $V^*$ in more generality.

Let the *valence* $d(v)$ of a vertex in a diagram $M$ be the number of edges incident to $v$, counting edges twice, that have both boundary vertices at $v$. In all of the following, when we speak of an *inner or interior vertex* of a diagram $M$ we mean an interior vertex of valence $d(v) \geq 3$, i.e. we delete all interior vertices of valence 2 in the sense of combining the adjacent edges of such a vertex to a single edge. (There are no interior vertices of valence 1 in a diagram, since relators are cyclically reduced) Under this assumption we define the degree $d(D)$ of a 2-cell $D$ of the diagram to be the number of edges in the boundary path of $D$, counting with multiplicity edges that are traversed twice (in opposite direction).

We now define the small cancellation conditions $W^*$ and $W$.

A diagram $M$ is said to be of *type* $W^*$ if it satisfies the following properties:

(a) Every 2-cell has degree greater than or equal to 3. (This corresponds essentially to the condition $C(3)$ of small cancellation theory)

(b) When every 2-cell $D$ of $M$ of degree $d$ is given the angles of a regular Euclidean $d$-gon, i.e. every corner of $D$ obtains angle $(1 - 2/d)\pi$, then at every inner vertex the sum of angles is greater or equal than $2\pi$.

A presentation $P$ is said to be of *type* $W^*$ if every vertex reduced diagram over $P$ is of type $W^*$.

We can express the condition $W^*$ combinatorially: $P$ is of *type* $W^*$ if for every vertex reduced diagram $M$ over $P$ and every inner vertex $v$ of $M$ the valence $d(v)$ and the degrees $d_i$ ($i = 1, ..., d(v)$) of the 2-cells incident with $v$ satisfy one of the following four conditions, where the cases $(i), (ii), (iii)$ are to be understood as follows: for each vertex of valence $d(v) = 3, 4, 5$ the degrees $d_i$ of adjacent 2-cells in a suitable order must be greater than or equal to the corresponding numbers in a single column of the table for $d_i$. (The degree of a 2-cell that has several corners at $v$ occurs with multiplicity):

(i) $d(v) = 3$,

| $d_1$ | 3  | 3  | 3  | 3  | 3  | 3  | 4  | 4  | 4  | 4  | 5  | 5  | 5  | 6 |
|-------|----|----|----|----|----|----|----|----|----|----|----|----|----|---|
| $d_2$ | 7  | 8  | 9  | 10 | 11 | 12 | 5  | 6  | 7  | 8  | 5  | 6  | 7  | 6 |
| $d_3$ | 42 | 24 | 18 | 15 | 14 | 12 | 20 | 12 | 10 | 8  | 10 | 8  | 7  | 6 |

(ii) $d(v) = 4$,

| $d_1$ | 3  | 3 | 3 | 3 | 3 | 4 |
|-------|----|---|---|---|---|---|
| $d_2$ | 3  | 3 | 3 | 4 | 4 | 4 |
| $d_3$ | 4  | 5 | 6 | 4 | 5 | 4 |
| $d_4$ | 12 | 8 | 6 | 6 | 5 | 4 |

(iii) $d(v) = 5$,

| $d_1$ | 3 | 3 |
|-------|---|---|
| $d_2$ | 3 | 3 |
| $d_3$ | 3 | 3 |
| $d_4$ | 3 | 4 |
| $d_5$ | 6 | 4 |

(iv) $d(v) \geq 6$,   $d_i \geq 3$ for $i = 1, ..., d(v)$

Note that, contrary to the standard non-metric small cancellation conditions (which are included as special cases), arbitrary combinations of different cases or subcases are permitted for the different vertices of the same diagram.

**Theorem 1** *Let $P$ be a finite presentation of type $W^*$. If $M$ is any vertex reduced disc-diagram over $P$ with $f$ 2-cells, then*

$$f \leq \frac{l^2(\delta M)}{\sqrt{3}\pi} \tag{1}$$

This theorem gives a quadratric isoperimetric inequality and, hence, a solution to the word problem.

In the statement of the following theorem we use the term "piece" as it is understood in small cancellation theory, i.e., in the context of diagrams over a presentation, a *piece* is a word that occurs as the label of an interior edge of some reduced diagram.

**Theorem 2** *Let $P = <x_1, \ldots, x_n \mid R_1, \ldots, R_m>$ be a finite presentation of type $W^*$ of the group $G$. Let $F$ be the free group on the generators and $p$ the maximal word length of pieces for the presentation $P$. If $u, v \in F$ are conjugate in $G$, then there exists a word $w \in F$ such that $u = wvw^{-1}$ and*

$$|w| \leq N + 2max\{|u|, |v|\}, \tag{2}$$

*where $N$ is the number of words in the alphabet $\{x_i^{\pm 1}\}$ of length $\leq 3 \cdot p \cdot max\{|u|, |v|\}$.*

This theorem gives a solution of the conjugacy problem by reducing it to the word problem.

One can dualize the definition of $W^*$ as follows:

A diagram is said to be of *type $W$* if it satisfies the following:

> If at every inner vertex of $M$ of valence $k$ (which is always $\geq 3$) every corner is given the angle $(2\pi)/k$, then the sum of the angles of every inner 2-cell $D$ is that of an Euclidean or hyperbolic polygon, i.e. less than or equal to $(d(D) - 2)\pi$. (An inner 2-cell is a closed 2-cell that is contained in the interior of $M$)

A finite presentation $P$ is said to be of *type $W$* if every <u>reduced</u> diagram over $P$ is of type $W$.

Combinatorially this is characterized by the same tables as above (in the case $W^*$) with valences of vertices and degrees of 2-cells interchanged, i.e. for every inner 2-cell $D$ of degree $d(D) = 3, 4, 5$, the valences $d_i$ of the vertices at the corners of $D$, in a suitable order, must be greater or equal to the numbers in one column of the corresponding table. If a vertex is incident to several corners of the same 2-cell D its degree occurs with multiplicity.

As for $W^*$-presentations, the word problem and conjugacy problem for groups which have presentations of type $W$ are solvable.

**Theorem 3** *Let $P$ be a finite presentation of type $W$, such that each generator occurs at least twice in the set $\{s_1, ..., s_m\}$ of the roots of relators (meaning precisely, that the number of occurences of letters from the set $\{x_i, x_i^{-1}\}$ in the collection of words is $\geq 2$). If $M$ is a reduced disk diagram over $P$ with $f$ 2-cells, then*

$$f \leq \frac{2(20^2)l^2(\delta M)}{\sqrt{3}\pi} \quad (3)$$

**Theorem 4** *Let $P = <x_1, \ldots, x_n \mid R_1, \ldots, R_m>$ be a finite presentation of type $W$ of the group $G$, where each generator occurs at least twice in the set $\{s_1, ...s_m\}$ of the roots of relators. Let $F$ be the free group on the generators and $r$ the maximal length of the relators of $P$. If $u, v \in F$ are conjugate in $G$, then there exists a word $w \in F$ such that $u = wvw^{-1}$ and*

$$|w| \leq N + 2max\{|u|, |v|\}, \quad (4)$$

*where $N$ is the number of words in the alphabet $\{x_i^{\pm 1}\}$ of length less than or equal to $31 \cdot r \cdot (|u| + |v|)$.*

**Remark:** The technical hypothesis in theorems 3 and 4 turns out to be not restrictive. Let $P$ be a presentation of type $W$ and suppose a generator, say $x_k$, occurs only once in the set of roots of the relators, i.e. $x_k$ does not occur in $R_j$ for $j \neq l$ and $R_l = (x_k^\epsilon w)^r$ and $x_k$ does not occur in $R_j$. Then, by a basis transformation, we see that $G = G' * \mathbb{Z}_r$ ($\mathbb{Z}_r$ is trivial if $r = 1$), where $G'$ has the presentation $P' = <x_i$ $(i = 1, \ldots, m, i \neq k) \mid R_j$ $(j = 1, \ldots, m, j \neq l) >$, which, as a subpresentation of $P$, is also of type $W$. But then the solvability of the word and conjugacy problems for $G$ follows from that for $G'$.

There is a similar restriction implicit in the definition of $W^*$ where we use the condition $C(3)$ "positively" (i.e. every relator must be a product of at least 3 pieces). In classical small cancellation theory it is defined negatively (no relator is a product of less than 3 pieces) including the possibility that it may not be a product of pieces at all, which happens exactly when a generator occurs only once in the set of roots of the relators. This is the only case where $W^*$ is not more general than the classical non-metric small cancellation conditions. The above discussion shows that this exception is irrelevant.

In [4] it is shown that presentations of type $W$ satisfy the cycletest and are therefore "combinatorially aspherical". This is also true by a quite similar proof for presentations of type $W^*$.

The following examples show that, in many cases, it may be quite easy to test if a presentation is of type $W$ or $W^*$.

**Example 1** *Let $P_n = <x, y, z \mid z^n = y, yx = xy>$ for $n \geq 2$. We will show that $P_n$ satisfies $W^*$ but not $W$.*

To prove $P$ is of type $W^*$ we must analyze the local situations around inner vertices of valence $\leq 6$ that can occur in a reduced diagram $f: M \to K_P$. The link of an inner vertex maps under $\pi \circ f$ onto a circuit $\gamma$ in the stargraph $S_P$. The length of $\gamma$ gives the valence of $v$, and $\gamma$ is reduced, i.e. has no backtracking, since $M$ is reduced.

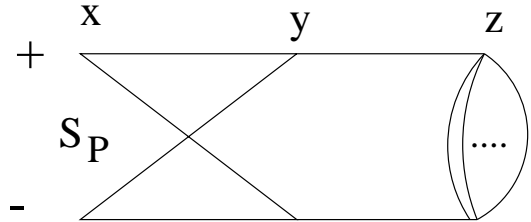

Figure 3: stargraph of $P$

The drawing of $S_P$ in figure 3 shows that there are two reduced circuits of length 5, one of length 4, and none of length 3. Figure 4 (b) shows the local situation around a vertex of valence 4. It is easy to see that the commutator relator can not break into pieces of length $> 1$. Therefore, the corresponding 2-cells have degree 4 and case (ii) of the definition of $W^*$ holds. Figure 4 (a) shows the local situation around a vertex corresponding to the circuit of length 5 that passes the vertices $+x, +y, +z, -z, -y$ in $S_P$ in that order. The edges of the 2-cells are already combined to pieces of maximal length. Therefore, the degrees of the adjacent 2-cells are 4, 4, and $\geq 3$ for the remaining three 2-cells, which is the second subcase of (iii) in the combinatorial definition of $W^*$. The other reduced circuit of length 5 yields essentially the same situation as in figure 4 (a), except that the edges labeled $x$ are oriented downwards.

It is clear that $P$ does not satisfy any of the standard small cancellation conditions. It is also not of type $W$. This can be seen by extending figure 4 (a) such that $Q$ and $R$ become inner vertices of a vertex reduced diagram, each having as link a circuit of length 5. Then the 2-cell $c$ of degree 3 has, at all its three corners, vertices of valence 5, contradicting the combinatorial definition of $W$.

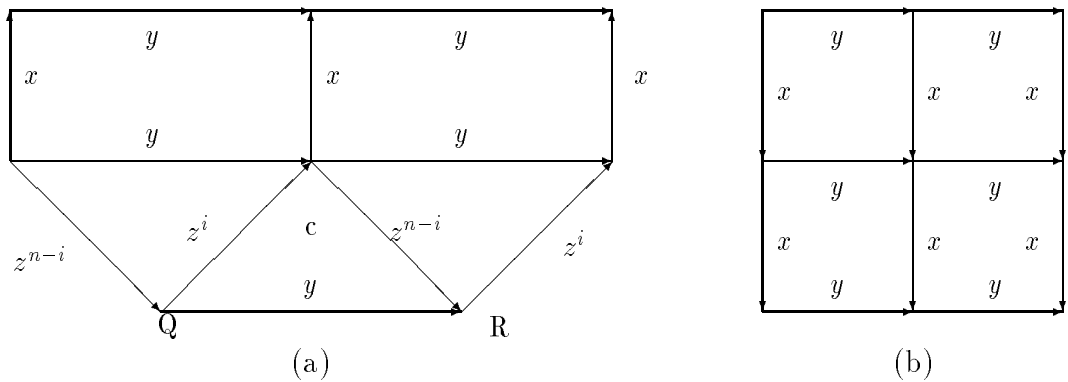

Figure 4: neighborhood of an inner vertex

**Example 2** *Let $P = <x, y \mid y^2x = xy^2>$. Juhász shows in [6] that this presentation satisfies $W$. $P$ is not of type $W^*$*

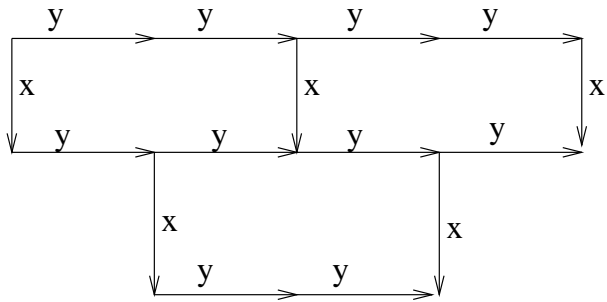

Figure 5: An inner vertex of valence 3

Figure 5 shows a neighborhood of an inner vertex of valence 3 in a reduced diagram over $P$. Since the only possible pieces for this presentation are $y^2$ and $y^{-2}$, the degrees of the adjacent 2-cells are 5 or 6. If the combinatorial definition of $W^*$ were true, then the degrees of all three 2-cells would have to be greater than or equal to 6. It is easy to extend Figure 5 to a larger diagram, where the three 2-cells have degree 5.

The group defined by $P$ is the same as the group defined by $P_2$ of the last example. This shows that $W$ and $W^*$ are attributes of the specific presentation and not of the group, which was to be expected.

## 2  Proof of Theorems 1 and 3

A combinatorial 2-complex is called *piecewise Euclidean* (PE) if each of its closed 2-cells has the metric of a convex polygon in the Euclidian plane and each of its closed 1-cells has the metric of a straight Euclidean line segment and these metrics agree on the overlaps. Let $M$ be a finite PE 2-disk, i.e. a PE 2-complex that is homeomorphic to a 2-dimensional disk. The given metric on each 2-cell permits to measure angles in the corners of the 2-cells of $M$. Let $v$ be an inner vertex of $M$ and let $g(v)$ be the sum of the angles that occur around $v$. Define the curvature at $v$ to be $\kappa(v) := 2\pi - g(v)$ and the *curvature of $M$* to be $\kappa(M) := \max \kappa(v)$, where the maximum is taken over all inner vertices $v \in M$ (set $\kappa(M) = 0$ if $M$ has no inner vertices). Let $l_M$ be the metric length of the boundary of $M$. The following theorem is due to Aleksandrov and Reshetnyak [10]:

**Theorem 5** *Let $M$ be a PE disk with $\kappa(M) \leq 0$, then*

$$Area(M) \leq \frac{l_M^2}{4\pi} \tag{5}$$

Here $Area(M)$ is the sum of the areas of the 2-cells of $M$. This result is certainly plausible. It gives as an upper estimate the area of a euclidean disk of the given boundary length. The idea of using the Theorem of Reshetnyak for proving a quadratic isoperimetric inequality for presentations of type $W$ is due to Gersten [3].

For the proofs of theorems 1 and 3 we may restrict ourselves to the case when the van Kampen diagram $M$ is a disk. The general case (for arbitrary v. Kampen diagrams) follows by induction from the following argument: If $M$ consists of two disks $M_1, M_2$ which are connected by some path or just by a common vertex and if theorems 1 and 3 are true for disks, i.e. $f_i \leq cl^2(\delta M_i)$ where $f_i$ is the number of 2-cells of $M_i$ ($i = 1, 2$) and $c$ is the constant of formula 1 or 3 respectively, then

$$f = f_1 + f_2 \leq c(l^2(\delta M_1) + l^2(\delta M_2)) \leq cl^2(\delta M).$$

**Proof of Theorem 1:** Assume $M$ is a vertex reduced disk-diagram over $P$ of type $W^*$ that is homeomorphic to a disk. We provide $M$ with the structure of a PE-2-complex (PE-structure) by realizing each closed 2-cell of $M$ of degree $d$ as a regular Euclidean $d$-gon with edges of length 1. (Recall that interior vertices of valence 2 are deleted) If the 2-cell $D$ has identifications on its boundary we must use a suitable stellar subdivision of a regular Euclidean $d$-gon for realizing $D$, so that any two edges or two vertices that get identified belong to different 2-cells of the subdivision. This is necessary since a closed 2-cell of a PE- structure can not have identifications on its boundary. By abuse of notation we call this PE-2-complex also $M$. The assumption that the diagram $M$ is of type $W^*$, guarantees that the sum of angles around each original interior vertex of $M$ is $\geq 2\pi$. Since the sum of angles at those interior vertices created by the subdivision is equal to $2\pi$, the PE-disk M has non-positive curvature and theorem 5 applies, yielding that $Area(M) \leq l^2(\delta M)/(4\pi)$. The smallest area of all regular n-gons of edge-length 1 is that of an equilateral triangle which has area $\sqrt{3}/4$. Therefore, if $f$ is the number of original 2-cells in $M$, $Area(M) \geq \sqrt{3}/4 f$ and the quadratic isoperimetric inequality (1) in theorem 1 follows from the fact that the length $l(\delta M)$ of the boundary word is equal to the metric length $l_M$ of the boundary of the PE-disk $M$. $\nabla$

In the proofs of theorems 3 and 4 we will make use of the following Curvature Lemma:

Let $M$ be a surface (with or without boundary) with a combinatorial cell decomposition, and let $g$ be an arbitrary *weight function* for $M$, i.e. a function from the set of corners of $M$ to the reals (One may think of the weight $g(\gamma)$ of a corner $\gamma$ as an angle which may assume arbitrary real values). We write $\gamma \prec D, \gamma \prec v$ to indicate that the corner $\gamma$ belongs to the 2-cell $D$, is incident with the vertex $v$, respectively. Let $g(D) := \Sigma_{\gamma \prec D} g(\gamma)$, $g(v) := \Sigma_{\gamma \prec v} g(\gamma)$, and let $\chi(M)$ denote the Euler characteristic of $M$. For $v \in M°$ (the interior of $M$) let $\kappa(v) := 2\pi - g(v)$, for $v \in \delta M$ let $\kappa(v) := \pi - g(v)$, and for a 2-cell $D \in M$ let $\kappa(D) := g(D) - (d(D) - 2)\pi$. Think of $\kappa(v)$, $\kappa(D)$ as the local curvature of $M$ at $v$, in $D$, respectively. Note that the local curvature at an interior vertex or in a 2-cell is zero exactly when the sum of angles is as in the euclidean case.

**Lemma 6 (Curvature–Lemma)** *If $g$ is an arbitrary weight function for a combinatorial cell decomposition of a surface $M$, then*

$$\sum_{v \in M} \kappa(v) + \sum_{D \in M} \kappa(D) = 2\pi \chi(M) \tag{6}$$

**Proof:** The proof is straightforward. Let $V, E, F$ be the number of vertices, edges, faces, respectively, of $M$, and let $V°, \dot{V}$, and $\dot{E}$ denote the number of interior vertices, boundary vertices, and boundary edges (i.e. edges in $\delta M$), respectively. Note that, since $M$ is a surface, $\delta M$ consists of disjoint circles, hence $\dot{E} = \dot{V}$.

By definition of $\kappa$:

$$\begin{aligned}
\sum_{v \in M} \kappa(v) &= \sum_{v \in M°} (2\pi - g(v)) + \sum_{v \in \delta M} (\pi - g(v)) \\
&= 2\pi V° + \pi \dot{V} - \sum_{v \in M} g(v) \\
&= \left(2V - \dot{V}\right)\pi - \sum_{v \in M} g(v)
\end{aligned} \tag{7}$$

$$\begin{aligned}
\sum_{D \in M} \kappa(D) &= \sum_{D \in M} [g(D) - (d(D) - 2)\pi] \\
&= \sum_{D \in M} g(D) - \sum_{D \in M} d(D)\pi + 2\pi F \\
&= \sum_{D \in M} g(D) - (2E - \dot{E})\pi + 2\pi F
\end{aligned} \tag{8}$$

By adding (7) and (8) and using $\sum_{v \in M} g(v) = \sum_{D \in M} g(D) =$ sum of the weights of all corners, we obtain the desired curvature formula:

$$\begin{aligned}
\sum \kappa(v) + \sum \kappa(D) &= (2V - \dot{V})\pi - \sum g(v) + \sum g(D) - (2E - \dot{E})\pi + 2\pi F \\
&= 2(V - E + F)\pi = 2\chi(M)\pi
\end{aligned}$$

$\nabla$

**Proof of Theorem 3:** Assume $P$ satisfies the hypothesis of theorem 3, and assume $M$ is a reduced diagram over $P$ (and hence of type $W$), such that $M$ is homeomorphic to a disk.

Let $M^*$ be the dual of $M$. Then $M^*$ consists, in general, of a collection of disks and arcs, where the disks are connected by arcs or by common points such that the result is simply connected; in other words, $M^*$ looks like a general disk diagram. Also, by the assumption that $M$ is of type $W$, $M^*$ is of type $W^*$. Hence, the proof of theorem 1 gives

$$f^* \leq l^2(\delta M^*)/(\sqrt{3}\pi) \tag{9}$$

for the number of 2-cells $f^*$ of $M^*$. We complete the proof of a quadratic isoperimetric inequality for the diagram $M$ itself by showing that $l(\delta M^*)$ is bounded by a linear function in terms of $l(\delta M)$ and by transforming (9) into an estimate for $f$ :

The assumption that each generator occurs at least twice in the set $\{s_1, ..., s_m\}$ of the roots of the relators, implies that the vertices of the star graph $S_P$ have valence $\geq 2$. It is then an easy exercise to show that every reduced path $w$ in $S_P$ can be extended to a reduced circuit $z$ containing $w$ as a subpath. This allows us to extend the entire diagram $M$ to a reduced diagram $M'$ over $P$, that contains $M$ in its interior and is, as a reduced diagram over $P$, also of type $W$. This is the point in the proof which prevents us from using "vertex reduced" in the definition of type-$W$ presentations; vertex reduced would not be preserved under the extension of $M$ to $M'$.

Recall that the combinatorial definition of $W$ lists lower bounds for the valences $d_i$ of the vertices of interior 2-cells of degrees 3, 4, and 5. The maximum of these lower bounds is 42, corresponding to an angle of $2\pi/42 = \pi/21$. If we define the angle of every corner $\gamma$ at an interior vertex $v_i$ of $M'$ to be $g(\gamma) = \max\{\pi/21, 2\pi/d(v_i)\}$ (where $d(v_i)$ is the valence of $v_i$ in $M'$), the sum of angles at every interior vertex of $M'$ will be $\geq 2\pi$, and the condition $W$ guarantees that for every interior 2-cell $D$ of $M'$ the sum of angles of its corners satisfies $\sum_{\gamma \in D} g(\gamma) \leq (d(D) - 2)\pi$. This holds in particular for every 2-cell of $M$. Under these conditions the curvature lemma for $M$ (with angle function induced from $M'$) reduces to

$$\sum_{v \in \delta M} \kappa(v) \left( = \sum_{v \in \delta M} (\pi - g(v)) \right) \geq 2\pi \tag{10}$$

or simply:

$$\sum_{v \in \delta M} g(v) \leq (l(\delta M) - 2)\pi \tag{11}$$

where the left hand side of (11) is the sum of angles of corners in $M$ that are incident with (vertices in) $\delta M$. Using that the number of corners incident with $\delta M$ equals $l(\delta M^*) + l(\delta M)$ and the fact that the angle of every corner is $\geq \pi/21$, we get from (11) the following estimate for $l(\delta(M^*))$:

$$l(\delta M^*) \leq 20 l(\delta M) - 42 \tag{12}$$

Taking into account that $V^\circ = V - l(\delta M) = f^*$ ($=$ the number of 2-cells of $M^*$) and that the Euler characteristic of $M^*$ is 2, one can derive from (9) the following estimate for the number $f$ of 2-cells of $M$:

$$f \leq \frac{2l^2(\delta M^*)}{\sqrt{3}\pi} + l(\delta M) - 2$$

By substituting (12) for $l(\delta M^*)$ and simplifying we obtain the quadratic isoperimetric inequality for $M$:

$$f \leq \frac{2(20^2)l^2(\delta M)}{\sqrt{3}\pi}$$

$\nabla$

**Remark:**  Details of the above proof are worked out, in principle, in the proof of theorem 5.9 in [5]. The quadratic isoperimetric inequality there is much sharper since we use a more restricted combinatorial definition of $W^*$ (which corresponds to the subcases of the definition above where the lower bounds for valences $d_i$ of vertices do not exceed 6). There may exist a more efficient way to estimate $l(\delta M^*)$ in terms of $l(\delta M)$. However, since the solution of the word problem by a quadratic isoperimetric inequality is rather theoretical, the size of the coefficient in the inequality does not matter that much.

## 3   Non-Positively Curved PE 2-complexes

Bridson studies in [1] the geometry of metric cell complexes under certain mild restrictions, which are certainly satisfied if the complexes are finite. His "Main Theorem" establishes the equivalence of certain local and global characterizations of non-positive curvature for simply connected metric cell complexes. In this section we will list the basic definitions and the results of [1] that we need in the proof of theorems 2 and 4. We will quote them as they apply to the special case of finite PE 2-complexes, which is a much narrower context than Bridson uses.

In a metric space $X$, a *geodesic segment* or, brief, a *geodesic* is a continuous path $\alpha$ that can be parametrized such that $\alpha\colon [0,l] \to X$ is an isometry from the interval $[0,l]$ to its image in X with the induced metric. $X$ is called a *geodesic metric space* if each pair of points $x, y \in X$ can be joined by a geodesic (in particular, such a space must be connected). Sometimes geodesic metric spaces are called path-metric spaces, since the distance between any two points can be realized as the length of a shortest path between them. A geodesic metric space $X$ is said to be *convex* if, for any pair of geodesics $\alpha, \beta$ in $X$, that are parametrized proportional to arc length on the interval $[0, 1]$, and for any $t \in [0, 1]$, the following inequality holds:

$$d(\alpha(t), \beta(t)) \leq (t-1)d(\alpha(0), \beta(0)) + td(\alpha(1), \beta(1)).$$

In the following we assume $K$ is a finite connected PE 2-complex. A *PL path* in $K$ is a finite concatenation of straight line segments $[x_i, x_{i+1}]$ where each line segment $[x_i, x_{i+1}]$ is contained in a 2-cell or a 1-cell of $K$. The length $l(\alpha)$ is the sum of the lengths of the line segments. The distance between two points $x, y \in K$ is defined to be $d(x, y) := inf\{l(\alpha) : \alpha$ a PL path from $x$ to $y\}$. Bridson shows that the infimum in the above definition is, in fact, a minimum; i.e. the distance between two points can be realized by a PL (geodesic) path. Hence, $d$ is a metric, called the *intrinsic metric* of the PE complex, and $K$, with this metric, is a geodesic metric space. The intrinsic metric is a little trickier than one would expect on first sight; for example, the metric induced by $d$ on a cell of $K$ does, in general, not coincide with the Euclidean metric on the cell (by which the PE complex was defined); however, this can be rectified by a suitably fine subdivision of the PE 2-complex, i.e., locally, the two metrics agree.

Let the length of an edge of the Whitehead graph $W(K)$ of a PE 2-complex $K$ be the angle at the corresponding corner of a 2-cell in $K$. A PE 2-complex $K$ is said to satisfy the *link condition* if the length of every non- trivial reduced circuit in $W(K)$ is greater or equal than $2\pi$. (A reduced circuit is a closed edge path that, considered as a cyclic path, has no backtracking.) This link condition defined by Bridson is a local condition of non-positive curvature. For a diagram $M$ provided with a PE structure, the link condition is equivalent to the sum of angles around every inner vertex of $M$ to be $\geq 2\pi$. It also coincides with the definition of non-positive curvature for PE 2-disks which we gave in section 2, above. Therefore, we will use the term "non-positively curved" for a PE 2-complex satisfying the link condition of Bridson.

The following two results from Bridsons article will be used in the proof of theorems 2 and 4. The first is part of the Theorem from section 2. of [1], the second is Lemma 2.3 (from the same section):

1. For a connected and simply connected PE 2-complex $K$ the following statements are equivalent:

    (I) $K$ is non-positively curved, i.e. $K$ satisfies the link condition,

    (II) $K$ has unique geodesics, i.e. for every pair of points $x, y \in K$ there is only one geodesic PL path connecting $x$ to $y$,

    (III) $K$ is convex (as a geodesic metric space).

2. If the PE 2-complex $K$ has unique geodesics then they vary continuously with their endpoints. More precisely, let $x_i \to x$ and $y_i \to y$ and let $\alpha_i, \alpha$ be the geodesics from $x_i$ to $y_i$, $x$ to $y$, respectively, then

    $$\|\alpha - \alpha_i\| = sup\{d(\alpha_i(t), \alpha(t)) : t \in [0, 1]\} \to 0.$$

    (Here "$\to$" denotes convergence of a sequence of points in $K$ or of a sequence of real numbers)

# 4   Proof of Theorems 2 and 4

Let $u, v$ be words in $F$. $u$ and $v$ are conjugate in $G$ if and only if there exists an annular diagram $M$ over $P$ whose boundary paths, oriented parallelly, read for suitable choices of start points the words $u$ and $v$. A *parallel orientation* is an orientation such that the boundary paths are freely homotopic in $M$. We will always assume that the boundary paths of an annular diagram are oriented parallelly.

One direction of this statement follows from the fact that $u = wvw^{-1}$ in $G$ corresponds to the equation $u = wvw^{-1} \prod_{j=1}^{k} w_j R_{i_j}^{\epsilon_j} w_j^{-1}$ in $F$. The right hand side of this equation can be realized by a disk diagram whose boundary reads $u$ and which has one exceptional region. The boundary of this exceptional region reads, instead of a relation of $P$, the word $v$. (In general, such a diagram may also have spherical components, attached by one-point-union. W.l.o.g. one may delete such spherical components - which only reduces the number $k$ of factors in the $\prod$-term of the equation - and obtain a disk diagram) Cutting out the exceptional region from the diagram then results in an annular diagram with boundary paths $u$ and $v$. Conversely, given such a diagram, select an edge path $w$ from the start point of $u$ to the start point of $v$. Then the closed path $uwv^{-1}w^{-1}$ is nullhomotopic in $M$ and therefore the word $uwv^{-1}w^{-1}$ is trivial in $G$.

**Proof of Theorem 2:**
Let $M$ be a vertex reduced annular diagram over $P$ whose boundary paths, oriented parallelly, read the words $u, v$. If the two boundary paths meet in at least one vertex, then some cyclic permutations of $u$ and $v$ will be equal in $G$, i.e. subwords of $u$, $v$ (or their inverses) respectively. This implies $u = a^{-1}bvb^{-1}a$ in $G$ with $|a^{-1}b| \leq 2 \max|u|, |v|$, and hence the conclusion of Theorem 2 holds. Therefore, we may restrict ourselves in the proof of Theorem 2 to the case when the boundary paths of $M$ are disjoint (as subsets of $M$). In this case $M$ consists of an actual annulus that has, in general, "trees of disks" attached to its boundaries (as sketched in figure 6).

Let $M'$ be the annulus obtained from $M$ by cutting off the trees of disks. To obtain a word $w$ that conjugates $v$ to $u$ we proceed as follows. Choose a shortest edge path $w'$ in $M'$ whose endpoints lie on opposite boundary components, then connect the endpoints of $w'$ along the boundary paths of $M$ to the start points of $u$ and $v$. This yields a conjugating word $w = xw'y$ where $x$ and $y$ are subwords of $u$ and $v$ (or their inverses), respectively, i.e. $|w| \leq 2\max\{|u|, |v|\} + |w'|$. It remains to show that, upon modifying the diagram $M'$, if necessary, without changing its boundary, we can achieve $|w'| \leq N$ where $N$ is the number of reduced words in $F$ of length $\leq 3 \cdot p \cdot \max\{|u|, |v|\}$. This follows directly from the subsequent Proposition applied to the diagram $M'$, taking into account that the lengths of the boundary paths of $M'$ are less than or equal to the lengths of the corresponding boundary paths of $M$.

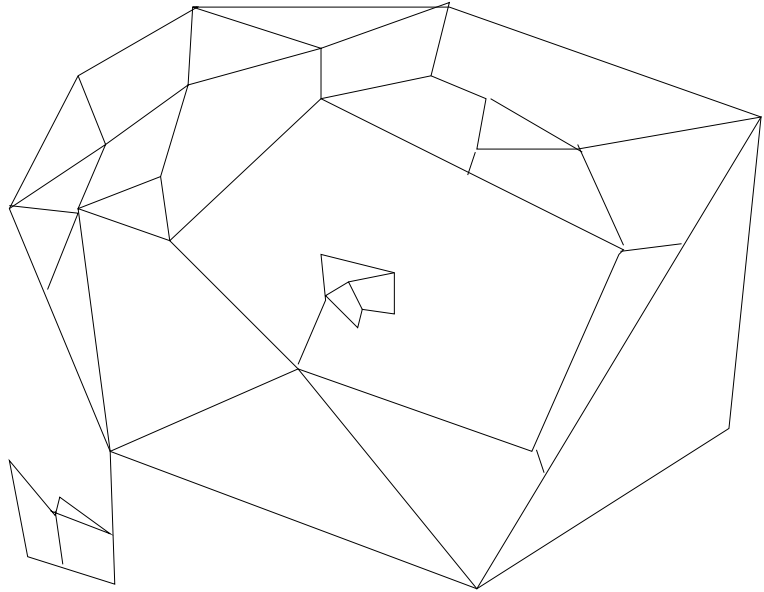

Figure 6: An annular diagram

**Proposition 7** *Let $M$ be an annular diagram of type $W^*$ over $P$, that is homeomorphic to an annulus. Let $l$ be the maximum of the lengths of its boundary paths and $N$ the number of words of length $\leq 3 \cdot p \cdot l$ over the alphabet $\{x_i^{\pm 1}\}$. Then there exists a word $w$ in $F$ of length $\leq N$ satisfying the equation $u = wvw^{-1}$ in $G$, where $u$ and $v$ are the words read from the boundary paths of $M$ for a suitable choice of startpoints.*

**Proof:** We provide the diagram $M$ with a PE-structure as in the proof of Theorem 1. I.e. every 2-cell $D$ of $M$ of degree n is represented by a regular Euclidean n-gon with edges of length 1. If $D$ has identifications on its boundary, we choose a stellar subdivision of the regular Euclidean n-gon so that the identifications of the boundary of $D$ do not contradict the Euclidean metric of the 2-cells of the subdivision. By the hypothesis that $P$ is of type $W^*$, the resulting PE-2-complex will be non-positively curved. By abuse of notation we call it also $M$, i.e. when we speak of the metric and the geodesics of $M$ we mean the metric and the geodesics of the PE-2-complex, which is actually defined on a subdivision of $M$. At the same time we will still think in terms of the original cell structure of $M$, and when we speak of vertices, edges, or 2-cells of $M$ we refer to the original (unsubdivided) cell structure. This does not create any serious problems. For example, geodesics intersect the original 2-cells in straight lines; more precisely, if one considers as a characteristic map for a 2-cell $D$ of $M$ a map $\phi$ whose domain is a regular Euclidean n-gon, then, for any geodesic $\gamma$ of $M$, the components of $\phi^{-1}(\gamma \cap D)$ are straight lines. This simply follows from the fact that, except at points that are vertices of the original cell structure, the PE-2-complex is locally isometric to the Euclidean plane or the Euclidean half plane.

In the following we will use $|\ |_m$ to denote the *metric length* of a path, i.e. the length of a PL-path or edge path in the metric of the PE 2- complex. $|\ |_w$ will denote the *word length* of an edge path, i.e. the length of the word which is the label of the edge path. For the boundary paths of $M$ the metric length coincides with the word length. It is clear that, in general, the word length of an edge path is less or equal than $p$ times its metric length, where $p$ is the maximal word length of a piece for the presentation $P$.

Now let $\bar{w}$ be a geodesic path in $M$ whose endpoints are vertices on the two boundary components of $M$, such that $\bar{w}$ realizes the shortest distance between the sets of vertices of the two boundary components. Let $Q$ be the startpoint and $R$ the endpoint of $\bar{w}$ and let $u, v$ be the closed boundary paths starting at $Q, R$ respectively. (The boundary words read along these paths will also be called $u$ and $v$.) We cut the annulus $M$ open along $\bar{w}$ to obtain $\bar{M}$. $\bar{M}$ inherits a PE- structure from $M$. (As mentioned above, $\bar{w}$ will cut through 2-cells of $M$ along straight lines creating pieces that are still convex Euclidean polygons) Since the link of an interior vertex of $\bar{M}$ is the same as the link of the corresponding vertex of $M$, $\bar{M}$ is a non-positively curved PE-disk. Hence, $\bar{M}$ has unique geodesic segments and is a convex geodesic metric space.

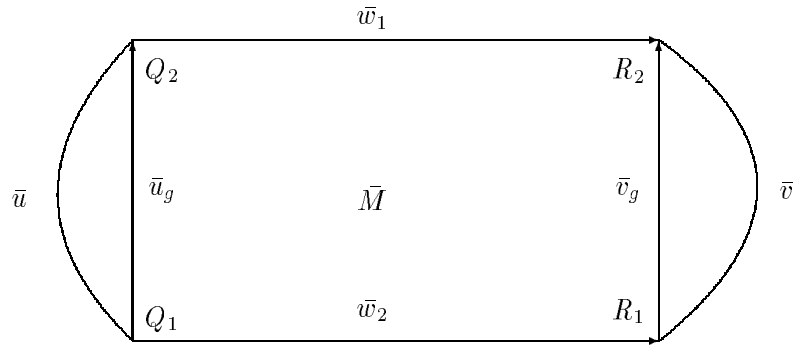

Figure 7: the diagram $M'$

For the following notation please compare figure 7. Let $Q_1$ and $Q_2$, $R_1$ and $R_2$, $\bar{w}_1$ and $\bar{w}_2$ be the pairs of vertices or paths in the boundary of $\bar{M}$ that originate from $Q$, $R$, $\bar{w}$, respectively, by the splitting of $M$ along $\bar{w}$, and let $\bar{u}$, $\bar{v}$ be the paths in the boundary of $\bar{M}$ that are created by cutting the closed boundary paths $u$, $v$ of $M$ at $Q$, $R$ respectively. We will also consider the (unique) geodesics $\bar{u}_g$ from $Q_1$ to $Q_2$ and $\bar{v}_g$ from $R_1$ to $R_2$ in $\bar{M}$. Together with $\bar{w}_1$ and $\bar{w}_2$ they form a geodesic rectangle in $\bar{M}$. If we parameterize each of the geodesics $\bar{w}_1$ and $\bar{w}_2$ proportional to arc length on the interval $[0, 1]$, we obtain by the convexity of $\bar{M}$

$$\begin{aligned}
d(\bar{w}_1(t), \bar{w}_2(t)) &\leq \max\{d(\bar{w}_1(0), \bar{w}_2(0)), d(\bar{w}_1(1), \bar{w}_2(1))\} \\
&= \max\{|\bar{u}_g|_m, |\bar{v}_g|_m\} \\
&\leq \max\{|u|, |v|\} = l,
\end{aligned}$$

where d denotes the distance in the geodesic metric space $\bar{M}$. For every $t \in [0,1]$ let $\bar{g}_t$ be the unique geodesic in $\bar{M}$ from $\bar{w}_1(t)$ to $\bar{w}_2(t)$ ($\bar{g}_0 = \bar{u}_g$, $\bar{g}_1 = \bar{v}_g$), and let $\bar{M}_1$ be the subspace of $\bar{M}$ bounded by the geodesic rectangle that is formed by $\bar{u}_g, \bar{v}_g, \bar{w}_1, \bar{w}_2$. Since in a simply connected PE-complex of non-positive curvature, geodesic segments vary continuously with their endpoints (see [1]), it follows that the one-parameter family of geodesics $\bar{g}_t$ covers all of $\bar{M}_1$, i.e. the map $\bar{g}: I \times I \to \bar{M}_1$ defined by $\bar{g}(t,s) = \bar{g}_t(s)$ is surjective.

Considering the annular diagram $M$ again, let $g_t$ be the family of closed PL-paths in $M$ that start and end at $\bar{w}(t)$ and are obtained from $\bar{g}_t$ by identifying the endpoints. It is easy to see that each path $g_t$ is shortest among the paths in $M$ that start and end at $\bar{w}(t)$ and are freely homotopic to $u$ (or $v$).

Now we replace the PL geodesic path $\bar{w}$ by a shortest edge path in $M$ from $Q$ to $R$, where "edge path" refers to the original cell structure of $M$. The label on $w$ reads a word in $F$ that conjugates $v$ to $u$ in $G$. Note that $w$ may intersect $u$ or $v$ in more than just the endpoints $Q$ or $R$ respectively and $w$ may not be homotopic relative its endpoints to $\bar{w}$. By the following Lemma 8 we obtain that $|w|_m \leq 2|\bar{w}|_m$.

**Lemma 8** *A locally geodesic PL-path $\alpha$ in $M$ that starts and ends at points $A$ and $B$ in $M^{(1)}$ (the 1-skeleton of $M$), respectively, is homotopic relative its endpoints to a PL-path $\beta$ in $M^{(1)}$ with the properties:*

1. $|\beta|_m \leq 2|\alpha|_m$,
2. *if $A$ and $B$ are in $M^{(0)}$ then $\beta$ is an edge path,*
3. *if $A$ or $B$ are not in $M^{(0)}$ then $\beta$ is an edge path except for the first or last segment.*

*We will call $\beta$ an edge path approximation of $\alpha$.*

**Proof:** The proof of Lemma 8 is elementary. Simply replace each straight line segment of $\alpha$ that intersects the interior of a 2-cell $D$ of $M$ by the shorter of the two paths on the boundary of $D$ that connect the same endpoints. Elementary geometry shows that the shorter path on the boundary of a regular Euclidean n-gon has at most two times the length of the secant through the interior of the polygon. After reducing any backtracking in the resulting path by a homotopy in the 1-skeleton we obtain a path $\beta$ that satisfies the conditions of Lemma 8. $\nabla$

**Lemma 9** *Let $w$ be the edge path constructed in the proof of Proposition 7 and let $V_0 = Q, V_1, V_2, ..., V_k = R$ be the sequence of vertices (including the interior vertices of valence 2) that are met by the path $w$. Then, for $i = 1, ..., k$, a shortest closed edge path $w_i$ in $M$, such that $w_i$ starts and ends at $V_i$ and is homotopic to $u$ (or $v$), will have metric length $|w_i|_m \leq 3 \cdot l$ and, hence, word length $|w_i|_w \leq 3 \cdot p \cdot l$.*

**Proof:** Recall that $\bar{M}_1$ is the part of $\bar{M}$ bounded by the geodesic rectangle with sides $\bar{u}_g, \bar{v}_g, \bar{w}_1, \bar{w}_2$. W.l.o.g. we can assume that $\bar{u}_g$ and $\bar{v}_g$ do not intersect; if they do, then it is easy to see that the edge path $w$, defined above, satisfies $|w|_m \leq 2\max\{|u|, |v|\}$ and hence the conclusion of Proposition 7. Let $M_1$ be the part of $M$ that corresponds to $\bar{M}_1$, i.e. $M_1$ is $\bar{M}_1$ with $\bar{w}_1$ and $\bar{w}_2$ identified. $M_1$ is an annulus with possibly one spike on either boundary. (The intersections $\bar{u} \cap \bar{w}_1$ and $\bar{u} \cap \bar{w}_2$, or $\bar{v} \cap \bar{w}_1$ and $\bar{v} \cap \bar{w}_2$ may consist of initial and terminal segments of $\bar{u}$ or $\bar{v}$, respectively, which create a spike)

If a vertex $V_i$ of the path $w$ is in $M_1$ then at least one of the closed PL-paths $g_t$, that were constructed above, will pass through $V_i$. Since $|g_t|_m \leq \max\{|u|, |v|\} = l$, a shortest closed PL-path in $M$ with basepoint $V_i$ that is homotopic to u will have length $\leq l$ and, by Lemma 8, a shortest edge path $w_i$ with the same properties will have length $|w_i|_m \leq 2l$.

If $V_i$ is not in $M_1$ we proceed as follows: Let $u_g$, $v_g$ be the closed PL-paths in $M$ that correspond to $\bar{u}_g$, $\bar{v}_g$ in $\bar{M}$; in other words, $u_g$ and $v_g$ are the boundary paths of $M_1$. Let $[t_a, t_b]$ be the maximal subinterval of $[0, 1]$ such that $w([t_a, t_b])$ is contained in $M_1$. Then $w(t_a)$ and $w(t_b)$ are points in $w \cap u_g$ and $w \cap v_g$ respectively. Since $w$ is a shortest edge path, Lemma 8 implies that the length of $w|_{[0, t_a]}$ is less or equal than two times the length of the shorter segment of $u_g$ connecting $w(0) = Q$ to $w(t_a)$. Hence $\left|w|_{[0,t_a]}\right|_m \leq |u_g|_m \leq |u|$. Similarly we obtain $\left|w|_{[t_b,1]}\right|_m \leq |v|$. Assume $V_i$ lies on the segment $w|_{[0,t_a]}$, i.e. $V_i = w(t_i)$ with $t_i < t_a$. Then $V_i$ can be connected to $Q = w(0)$ along the edge path $\lambda_i = (w|_{[0,t_i]})^{-1}$ and $|\lambda_i|_m \leq |u|$. Hence, $\lambda_i u \lambda_i^{-1}$ is a closed edge path homotopic to $u$ with basepoint $V_i$ that has metric length $\leq 3|u| \leq 3l$. If $V_i$ lies on $w|_{[t_b,1]}$ we choose $\lambda_i$ to be a final segment of $w$ connecting $V_i$ to $R = w(1)$, and the path $\lambda_i v \lambda_i^{-1}$ will have the desired properties and length estimate. $\nabla$

We are now able to complete the proof of Proposition 7 and thereby the proof of Theorem 2.

Consider the diagram M of Proposition 7 and the edge path $w$ in $M$ from $Q$ to $R$. Let $w_{QV_i}$, $w_{V_j R}$ be the segments of $w$ from $Q$ to $V_i$, $V_j$ to $R$, respectively, and let $w_i$ be the paths from Lemma 9. Any null homotopic edge path in a diagram represents a word that is equal to 1 $G$ and hence corresponds to a closed path in the Cayley graph. This applies to the closed edge paths $u w_{QV_i} w_i^{-1} w_{QV_i}^{-1}$ and $w_j w_{V_j R} v w_{V_j R}^{-1}$. We now consider the conjugacy relation $uwv^{-1}w^{-1}$ as a closed path based at 1 in the Cayley graph of $G$.

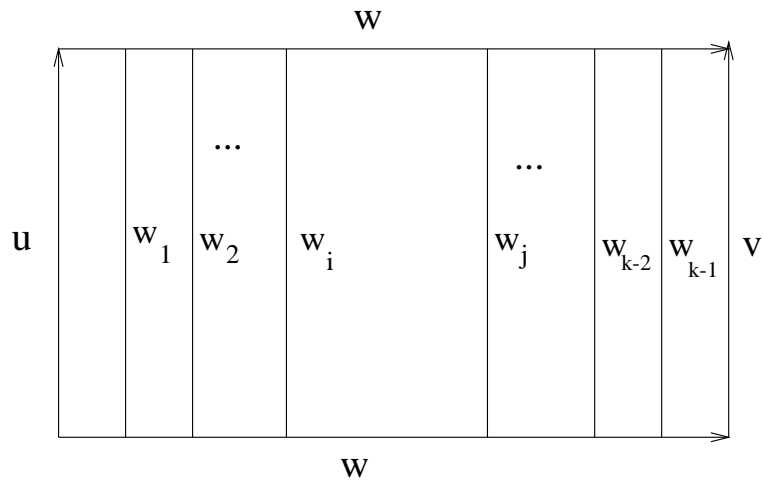

Figure 8: A part of the Cayley graph

We then have the situation depicted in figure 8, and we are able to apply a standard surgery trick: Let $N$ be the number of words of length $\leq 3pl$ over the alphabet $\{x_1^{\pm 1}, \ldots, x_n^{\pm 1}\}$. If $|w|_w > N$ then some of the "vertical words", say $w_i$ and $w_j$ will be equal, and we can do surgery on the diagram of figure 8 by cutting out the middle part and gluing the shaded left part to the shaded right part along $w_i = w_j$. This is done in the Cayley graph by translating the right part under left multiplication by $w_{V_i V_j}^{-1}$. (Here we use that left multiplication by any group element is an isometry of the Cayley graph that preserves edge labels.) This surgery shortens the conjugating word $w$ by the length of the word $w_{V_i V_j}$. Since the vertical paths in the resulting figure are those in the shaded parts of figure 8, we can iterate this procedure until the length of the horizontal (conjugating) word becomes less than $N$.    $\nabla$

**Proof of Theorem 4:**
We follow the basic idea of the proof of Theorem 2 where the key steps are Proposition 7 and Lemma 9. Given a presentation $P$ that satisfies the hypothesis of Theorem 4, let $M$ be a reduced annular diagram over $P$ that is homeomorphic to an annulus and let $M^*$ be the dual of $M$. Since $M$ is of type $W$, $M^*$ is of type $W^*$.

We first wish to estimate the length of the boundary $l(\delta M^*)$ in terms of $l(\delta M)$ (where the length is measured, as usual, as the number of edges in the boundary paths). We use the same trick as in the proof of Theorem 3 to do so: Extend the diagram $M$ to a diagram $M'$ that contains $M$ in its interior and is also reduced and hence of type $W$. If we give every corner $\gamma$ at every interior vertex of $M'$ the angle $g(\gamma) = \max\{\pi/21, 2\pi/d(v_i)\}$ then the sum of angles around every interior vertex is $\geq 2\pi$ and the sum of angles of every interior 2-cell $D$ of $M'$ and hence every 2-cell of $M$ is $\leq (d(D) - 2)\pi$. Since $\chi(M) = 0$ we obtain from the curvature lemma that the sum of angles of the corners of $\delta M$, i.e. the corners of $M$ that are incident with vertices in $\delta M$, is less than or equal to $l(\delta M)\pi$. Taking into account that the minimal angle of a corner is $\pi/21$, we get:

$$\#\text{of corners of } \delta M = l(\delta M^*) + l(\delta(M) \leq 21 l(\delta M) \qquad \text{or:} \qquad l(\delta M^*) \leq 20 l(\delta M).$$

As in the proof of Theorem 2 we may assume, w.l.o.g., that $M^*$ consists of an actual annulus $M^{*\prime}$ that has (in general) trees of disks attached to its boundaries. We make $M^{*\prime}$ into a PE 2-complex and construct an edge path $w^*$ connecting vertices of the opposite boundary components of $M^{*\prime}$ as in the proof of Proposition 7. Then the conclusion of Lemma 9 holds for $w^*$ in $M^{*\prime}$, i.e. for every vertex $V_i^*$ of $w^*$ there exists a closed edge path $w_i^*$ such that $w_i^*$ starts and ends at $V_i^*$, is homotopic to a boundary path of $M^{*\prime}$, and $|w_i^*|_m \leq 3l(\delta M^{*\prime}) \leq 3l(\delta M^*)$.

By duality, any edge path $\alpha^*$ in $M^*$ corresponds to a sequence $D_0 e_0 D_1 e_1 ... e_{k-1} D_k$ of 2-cells $D_i$ and 1-cells $e_i$ in $M$ such that $e_i \subset D_i \cap D_{i+1}$ (we call it briefly a *chain of 2-cells*). Clearly, for such a chain of 2-cells in $M$ belonging to an edge path $\alpha^*$ in $M^*$ and for any pair of vertices $Q$ in $D_0$ and $R$ in $D_k$ one can find an edge path $\alpha$ from $Q$ to $R$ in $M$ that is contained in the subcomplex $\bigcup_{i=1}^{k} D_i$ and has word-length $|\alpha| \leq r(|\alpha^*| + 1)/2$, where $r$ is the maximal length of the relators of $P$. We call $\alpha$ an *accompanying path* of $\alpha^*$. First we apply this to the path $w^*$ and obtain an accompanying edge path $w$ from $Q$ to $R$ in $M$, where we can choose $Q$ and $R$ to be vertices in opposite boundary components of $\delta M$. (This follows since $Q^*$ and $R^*$ are vertices in opposite boundary components of $\delta M^{*\prime} \subset \delta M^*$ and a vertex of $M^*$ belongs to a boundary component of $M^*$ if and only if the corresponding 2-cell of $M$ intersects the corresponding boundary component of $M$.) By deleting loops or backtracking we may assume that $w$ is a simple path, i.e. without double points. Call the vertices along $w$ : $V_i$. Then each $V_i$ belongs to a 2-cell $D_{j_i}$ of the chain corresponding to $w^*$ and we can find an accompanying path $w_i$ of $w_{j_i}^*$ with startpoint and endpoint $V_i$ that is freely homotopic to either boundary component of $M$ and has length

$$|w_i|_w \leq \frac{r(|w_{j_i}^*| + 1)}{2} \leq \frac{r(3l(\delta M^*) + 1)}{2} \leq 31 \cdot r \cdot l(\delta M).$$

Now we are ready to do the same surgery trick as in the proof of Proposition 7 and obtain the following result: If $u$, $v$ are the words read from the boundary of $M$ starting at $Q$, $R$, respectively, and $N$ is the number of words in the alphabet $\{x_1^{\pm 1}, \ldots, x_n^{\pm 1}\}$ of length $\leq 31 r l(\delta M)$, then there exists a word $w$ in $F$ of length $\leq N$ such that $u = wvw^{-1}$ holds in $G$.

Finally, if we have an arbitrary reduced annular diagram $M$ with given boundary words $u$ and $v$ and $M$ is not homeomorphic to an annulus, then, by the same argument as in the proof of Theorem 2, there exists a word of length $\leq 2\max\{|u|, |v|\} + N$ conjugating $u$ to $v$ in $G$. $\qquad \nabla$

# 5   The Conditions $V$ and $V^*$

We define generalized small cancelation conditions that imply a linear isoperimetric inequality and, hence, word-hyperbolicity for the group. The situation will be simpler since the conjugacy problem is obtained for free. We start with very general conditions on a presentation $P$ that imply a linear estimate of the number of vertices or the number of 2-cells of a reduced diagram in terms of the length of the boundary of the diagram. Let $g$ be a weight function that assigns real valued weights to the corners of a diagram $M$. Here we think of the weights, multiplied by $\pi$, as angles. Note that this differs slightly from section 2, where we had the factor $\pi$ incorporated into the weights. Using weights in diagrams to determine hyperbolicity of groups goes back to Gersten [2] and Pride [8] and the proofs exhibited here are very similar to those in [2] and [8]. For a vertex $v$ and a 2-cell $D$ let $g(v)$, $g(D)$ be defined as in section 2 and let $g(\delta M) := \Sigma_{v \in \delta M} g(v)$.

**Theorem 10** *Given a presentation $P$. If there exist real numbers $\epsilon > 0$ and $N \geq 0$ such that for every vertex reduced diagram $M$ over $P$ there is a weight function $g$ satisfying*

1. *$g(D) \leq d(D) - 2$ for every 2-cell $D$,*

2. *$g(v) \geq 2 + \epsilon$ for every inner vertex $v$, and*

3. *$g(\delta M) \geq -N\ l(\delta M)$,*

*then the number $V$ of vertices of $M$ is bounded by*

$$\begin{array}{rll} V \leq & \frac{2+N}{\epsilon} l(\delta M) & \text{if } \chi(M) \geq 0, \text{ and} \\ V \leq & \frac{4+N}{\epsilon} l(\delta M) & \text{if } \chi(M) < 0. \end{array}$$

**Proof:**   W.l.o.g. we assume $\epsilon < 1$. Let $V^\circ, E, F$, be the numbers of interior vertices, edges, faces, respectively, of $M$, and let $s = \Sigma_{D \in M} g(D) = \Sigma_{v \in M} g(v)$ be the sum of weights of all corners of $M$. From 1. and the fact that $\Sigma_{D \in M} d(D) = 2E - l(\delta M)$ (which is easy to verify) we obtain

$$s = \sum_{D \in M} g(D) \leq \sum_{D \in M} d(D) - 2F = 2E - l(\delta M) - 2F.$$

On the other hand, using 2. and 3. we get :

$$s = \sum_{v \in M} g(v) = \sum_{v \in M^\circ} g(v) + \sum_{v \in \delta M} g(v) \geq (2 + \epsilon) V^\circ - N \cdot l(\delta M),$$

where $M^\circ$ denotes the interior of $M$. Combining both inequalities and using $V^\circ \geq V - l(\delta M)$ gives:

$$\begin{array}{rcl} (2 + \epsilon)(V - l(\delta M)) - Nl(\delta M) & \leq & 2E - l(\delta M) - 2F \\ 2V - 2E + 2F + \epsilon V - (2 + \epsilon) l(\delta M) & \leq & -l(\delta M) + Nl(\delta M) \\ 2\chi(M) + \epsilon V & \leq & (1 + \epsilon + N) l(\delta M). \end{array}$$

When $\chi(M)$ is non-negative (i.e. when $M$ is a disk diagram or an annular diagram) we get

$$V \leq \frac{1+\epsilon+N}{\epsilon} l(\delta M) \leq \frac{2+N}{\epsilon} l(\delta M).$$

In general, we estimate $\chi(M)$ by $\chi(M) = 2 - \tau \geq 2 - l(\delta M)$ (where $\tau$ is the number of boundary paths of the diagram $M$) and get:

$$\begin{aligned} 2(2 - l(\delta M)) + \epsilon V &\leq (1+\epsilon+N)l(\delta M) \\ 4 + \epsilon V &\leq (3+\epsilon+N)l(\delta M) \\ V &\leq \frac{4+N}{\epsilon} l(\delta M) \end{aligned}$$

$\nabla$

**Theorem 11** *Given a presentation $P$. If there exist real numbers $\epsilon > 0$ and $N \geq 0$ such that for every reduced diagram $M$ over $P$ there is a weight function $g$ satisfying*

1. $g(D) \leq d(D) - 2 - \epsilon$ *for every 2-cell $D$,*
2. $g(v) \geq 2$ *for every inner vertex $v$, and*
3. $g(\delta M) \geq -N\ l(\delta M)$,

*then the number $F$ of 2-cells of $M$ is bounded by*

$$F \leq \frac{3+N}{\epsilon}\ l(\delta M).$$

The proof of Theorem 10 is very similar to the proof of Theorem 9. It is given as proof of Theorem 5.4 in [5] and originally due to Pride [8]. Note that the inequalities in the conclusion of Theorems 10 and 11 apply to arbitrary diagrams $M$, not just disk diagrams.

A strengthening of the conditions $W^*$ and $W$ leads to classes of presentations $V^*$ and $V$, defined below, which satisfy the hypotheses of Theorems 10 and 11, respectively. Groups with presentations of type $V^*$ or $V$ are word-hyperbolic since they satisfy a linear isoperimetric inequality (Theorem 12).

A diagram $M$ is said to be *of type $V^*$* if it is of type $W^*$ with the sum of angles around every interior vertex being strictly greater than $2\pi$. A presentation $P$ is *of type $V^*$* if every vertex reduced diagram $M$ over $P$ is of type $V^*$. The combinatorial definition of $V^*$ is given in the same way as the combinatorial definition of $W^*$ by the following tables (where $d_i$ are lower bounds for the degrees of 2-cells adjacent to an interior vertex $v$):

(i) $d(v) = 3$,

| $d_1$ | 3 | 3 | 3 | 3 | 3 | 3 | 4 | 4 | 4 | 4 | 5 | 5 | 5 | 6 |
|---|---|---|---|---|---|---|---|---|---|---|---|---|---|---|
| $d_2$ | 7 | 8 | 9 | 10 | 11 | 12 | 5 | 6 | 7 | 8 | 5 | 6 | 7 | 6 |
| $d_3$ | 43 | 25 | 19 | 16 | 14 | 13 | 21 | 13 | 10 | 9 | 11 | 8 | 7 | 7 |

(ii) $d(v) = 4$,

| $d_1$ | 3 | 3 | 3 | 3 | 3 | 4 |
|---|---|---|---|---|---|---|
| $d_2$ | 3 | 3 | 3 | 4 | 4 | 4 |
| $d_3$ | 4 | 5 | 6 | 4 | 5 | 4 |
| $d_4$ | 13 | 8 | 7 | 7 | 5 | 5 |

(iii) $d(v) = 5$,

| $d_1$ | 3 | 3 | 3 |
|---|---|---|---|
| $d_2$ | 3 | 3 | 3 |
| $d_3$ | 3 | 3 | 4 |
| $d_4$ | 3 | 4 | 4 |
| $d_5$ | 7 | 5 | 4 |

(iv) $d(v) = 6$,    $d_i \geq 3$ for $i = 1, ..., d(v)$ and $d_i > 3$ for at least one $i$

(v) $d(v) \geq 7$,    $d_i \geq 3$ for $i = 1, ..., d(v)$.

A diagram $M$ is said to be *of type V* if it is of type $W$ with the sum of angles of every interior 2-cell $D$ being strictly less than $d(D) - 2$. A presentation $P$ is *of type V* if every reduced diagram over $P$ is of type $V$. Combinatorially, type $V$ is characterized by the same tables as above where valences of vertices and degrees of 2-cells are interchanged and the vertices and 2-cells are both from the interior of the diagram.

**Theorem 12** *Assume $G$ is a group that has a finite presentation of type $V^*$ or of type $V$ and, in the latter case, assume that each generator occurs at least twice in the set of roots of the relators. Then $G$ satisfies a linear isoperimetric inequality and hence is word-hyperbolic.*

**Proof:** We consider first a presentation $P$ of type $V^*$. From the definition of $V^*$ it follows easily that $P$ satisfies the hypothesis of Theorem 10: Let $M$ be a vertex reduced disk diagram over $P$. Define the weight of a corner of $M$ to be the angle given by the definition of $V^*$, divided by $\pi$. Then $g(D) = d(D) - 2$ for every 2-cell $D$ of $M$ and $g(v) > 2$ for every interior vertex $v$ of $M$ and, since this latter requirement is characterized combinatorially by the finitely many cases listed above, it is clear that there must exist an $\epsilon > 0$ such that $g(v) \geq 2 + \epsilon$ for interior vertices. More precisely, we can choose $\epsilon$ to be $1/903$. This minimal "weight excess" occurs when $v$ has valence 3 and the degrees $d_i$ of adjacent 2-cells are 3, 7, and 43 corresponding to weights $1/3$, $5/7$, and $41/43$ which add up to $2 + 1/903$. Note also that all weights are positive, hence we can choose $N$ to be 0.

We get from Theorem 10 the estimate $V \leq \frac{2}{\epsilon} l(\delta M)$ for the number $V$ of vertices of $M$, which transforms into a linear estimate for the number $F$ of 2-cells of $M$ as follows: By the definition of $V^*$ we know that every 2-cell of $M$ has degree greater than or equal to 3, hence

$$3F \leq \sum_{D \in M} d(D) = 2E - l(\delta M)$$

$$3F - 2E + 2V \leq 2V - l(\delta M)$$
$$F + 2\chi(M) \leq \left[\frac{4}{\epsilon} - 1\right] l(\delta M)$$

The last inequality uses the estimate $V \leq \frac{2}{\epsilon} l(\delta M)$ and gives the desired linear isoperimetric inequality for $F$ since the Euler characteristic for a disk diagram is 1.

To see that a presentation $P$ of type $V$ satisfies the hypothesis of Theorem 11 and hence the linear isoperimetric inequality which is the conclusion of Theorem 11, requires a little extra work. (In the definition of the class $V$ only corners at interior vertices are given angles, whereas for Theorem 11 every corner of the diagram needs to be given a weight) This can be solved by the same extension trick that was used in the proofs of Theorems 3 and 4: Given a reduced disk diagram $M$ over $P$, cut $M$ at vertices where it fails to be a 2-dimensional manifold (with boundary), thereby separating the diagram into pieces which are either 2-dimensional disks or 1-dimensional arcs. Call the pieces that are disks $M_i$. Now extend each diagram $M_i$ along its boundary to a reduced diagram $M'_i$, as in the proof of Theorem 3, such that all vertices and all 2-cells of $M_i$ lie in the interior of the extended diagram $M'_i$. Since $P$ is of type $V$, $M'_i$, as a reduced diagram over $P$, is of type $V$. Therefore, if we give any corner at an interior vertex $v$ of $M'_i$ weight $\frac{2}{d(v)}$, the sum of weights in every interior 2-cell $D$ of $M'_i$ will be less than $d(D) - 2$. Now, considering the subdiagrams $M_i$ (with the induced weights from $M'_i$) again as parts of the original diagram $M$, we have obtained a positive weight function for the corners of $M$ satisfying the hypothesis of Theorem 11. (That the sum of weights of every 2-cell is actually less or equal than $d(D) - 2 - \epsilon$ for $\epsilon = 1/903$ follows, in the same manner as above, from the combinatorial characterization of "type $V$" by finitely many cases. We conclude that a presentation $P$ of type $V$ satisfies the linear isoperimetric inequality from the conclusion of Theorem 11 with $N = 0$ and $\epsilon = 1/903$. $\nabla$

**Remarks:**

1. The proof above applies with minor modifications also to diagrams with more than one boundary path, yielding a linear estimate for the number of 2-cells in terms of the length of the boundary.

2. We wish to emphasize that Theorem 10 and 11 can be used to prove that more general classes of presentations define hyperbolic groups. The proofs concerning the small cancellation classes $W^*$ and $W$, which use metric 2-complexes of non-positive curvature as the main tool, depend strongly on the 2-cells in a diagram or in the dual of a diagram being realized as regular polygons. In the negatively curved case the proofs are easier, using just Euler characteristic counts, and do not depend on regular distributions of angles in 2-cells or around vertices. Hence there are other classes of presentations allowing "non regular" weight functions on diagrams to which the proofs above apply. Such classes are the presentations satisfying the hyperbolic weight test as defined by Gersten [2] and extended by Pride [8] and, more generally, the presentations satisfying the hyperbolic cycle test [4]. The latter include the presentations of type $V$ and, using a dual form of the hyperbolic cycle test, the presentations of type $V^*$.

Stephan Rosebrock        Günther Huck
Hauptstr. 38             Dept. of Math.
65779 Kelkheim           Northern Arizona University
Germany                  Flagstaff AZ 86011
                         USA